\documentclass[12pt]{amsart}
\usepackage{amsmath,amscd,amssymb,amsfonts,graphics}
\newcommand{\h}{\hbox}
\newcommand{\q}{\quad}
\newcommand{\nin}{\noindent}
\newcommand{\bs}{\par\bigskip}
\newcommand{\ms}{\par\medskip}
\newcommand{\sk}{\par\smallskip}
\newcommand{\bsn}{\par\bigskip\noindent}
\newcommand{\msn}{\par\medskip\noindent}
\newcommand{\skn}{\par\smallskip\noindent}
\newcommand{\ges}{\geqslant}
\newcommand{\les}{\leqslant}
\newcommand{\1}{\hskip1pt}
\newcommand{\mcap}{\hbox{$\bigcap$}}
\newcommand{\mcup}{\hbox{$\bigcup$}}

\newcommand{\msum}{\hbox{$\sum$}}

\newcommand{\mprod}{\hbox{$\prod$}}
\newcommand{\A}{{\mathcal A}}

\newcommand{\G}{{\mathcal G}}
\newcommand{\I}{{\mathcal I}}
\newcommand{\Hc}{{\mathcal H}}
\newcommand{\Lc}{{\mathcal L}}
\newcommand{\OO}{{\mathcal O}}
\newcommand{\Sc}{{\mathcal S}}

\newcommand{\PP}{{\mathbb P}}
\newcommand{\Q}{{\mathbb Q}}
\newcommand{\C}{{\mathbb C}}

\newcommand{\R}{{\mathbb R}}
\newcommand{\RR}{{\mathbf R}}

\newcommand{\Z}{{\mathbb Z}}

\newcommand{\Gr}{{\rm Gr}}

\newcommand{\Et}{\widetilde{E}}
\newcommand{\Yt}{\widetilde{Y}}
\newcommand{\Ft}{\widetilde{\mathcal F}}
\newcommand{\jt}{{}\,\widetilde{\!j}{}}
\newcommand{\Ff}{F_{\!f}}
\newcommand{\al}{\alpha}
\newcommand{\la}{\lambda}

\newcommand{\om}{\omega}

\newcommand{\DR}{{\rm DR}}
\newcommand{\ord}{{\rm ord}}

\newcommand{\mult}{{\rm mult}}

\newcommand{\ddd}{{\rm d}}
\newcommand{\bl}{\bigl}
\newcommand{\br}{\bigr}
\newcommand{\pl}{\1{+}\1}
\newcommand{\mi}{\1{-}\1}
\newcommand{\eq}{\,{=}\,}

\newcommand{\less}{\,{\leqslant}\,}
\newcommand{\gess}{\,{\geqslant}\,}

\newcommand{\ssb}{\raise.15ex\h{${\scriptscriptstyle\bullet}$}}
\newcommand{\ssc}{\,\raise.15ex\h{${\scriptstyle\circ}$}\,}

\newcommand{\into}{\hookrightarrow}
\newcommand{\simto}{\,\,\rlap{\hskip1.5mm\raise1.4mm\hbox{$\sim$}}\hbox{$\longrightarrow$}\,\,}
\newcommand{\plim}{\rlap{\raise-5.5pt\h{$\,\leftarrow$}}{\rm lim}}
\begin{document}
\title[Efficiency and complexity of arrangements]
{Efficiency and complexity of hyperplane arrangements}
\author[M. Saito]{Morihiko Saito}
\address{M. Saito : RIMS Kyoto University, Kyoto 606-8502 Japan}
\email{msaito@kurims.kyoto-u.ac.jp}
\begin{abstract} For a projective hyperplane arrangement, we study sufficient conditions in terms of combinatorial data for ESV-calculability of the monodromy eigenspaces of the first Milnor fiber cohomology for eigenvalues of order $m>1$. This can be reduced to the line arrangement case by Artin's vanishing theorem. These sufficient conditions are often unsatisfied if efficiency or complexity of the combinatorics of arrangement is high. In order to measure these, we introduce the notions of $m$-efficiency and $m$-complexity for $m\gess 3$. The former is defined to be the number of points with multiplicity divisible by $m$ lying on one line in average. In many cases, one of the above sufficient conditions is satisfied if it is at most 2, although there are certain exceptional cases, especially when $m\eq 3$. The $m$-complexity is defined to be the maximal number of edges containing one vertex of the associated $m$-graph. We can show that one of the sufficient condition holds if it is at most $(m{+}1)/2$.
\end{abstract}
\maketitle
\ms\centerline{\bf Introduction}
\bsn
Let $X\subset\PP^{n-1}$ be a reduced projective hyperplane arrangement with $f$ a defining polynomial of $n$ variables. To calculate a monodromy eigenspace of its first Milnor fiber cohomology $H^1(\Ff,\C)_{\la}$ for eigenvalue $\la$ of order $m>1$, we may {\it assume\1} $n\eq 3$ (where $X$ is denoted by $L$). Indeed, the computation is reduced to this case using an iterated general hyperplane cut together with Artin's vanishing theorem \cite{BBD}. Note that the dimension of the second Milnor fiber cohomology is then determined from the first together with the Euler number of $\PP^2\setminus L$ using (1.2.1) below.
\sk
We may also assume $\tfrac{d}{m}\in\Z$ with $m:=\ord\,\la$, since the eigenspace $H^j(\Ff,\C)_{\la}$ vanishes unless $\tfrac{d}{m}\in\Z$. These eigenspaces can be calculated by using the corresponding {\it Aomoto complex\1} if some condition coming from \cite{ESV} is satisfied, see \cite{BDS}, \cite{BSY}, \cite{bha}, etc. Combining this with the theory of cyclotomic polynomials, we can get the following (see (1.2) below).
\msn
{\bf Theorem~1.} {\it Let $L\subset\PP^2$ be a reduced line arrangement defined by a homogeneous polynomial $f$ of $3$ variables with degree $d$. Then the eigenspace $H^1(\Ff,\C)_{\la}$ with eigenvalue $\la$ of order $m>1$ is ESV-calculable, that is, it can be calculated by the corresponding Aomoto complex, if there is a subset $J\subset\{1,\dots,d\}$ with $|J|=\tfrac{d}{m}$, and one of the following two conditions is satisfied\,$:$}
$$L^{[\ges 3]}\cap L^{[(m)]}\cap L_J^{[k]}\subset L^{[\ges km]}\q(\forall\,k\ges 2),
\leqno{\rm (a)}$$
\vskip-6mm
$$L^{[\ges 3]}\cap L^{[km]}\subset L^{[\ges k]}_J\q(\forall\,k\ges 1).
\leqno{\rm (b)}$$
\sk
Here $L_J:=\mcup_{i\in J}\,L_i$ with $L_i\,\,(i\in\{1,\dots,d\})$ the irreducible components of $L$, and
$$\aligned L^{[k]}&:=\{P\in L\mid\mult_PL=k\},\q L^{[\ges k]}:=\{P\in L\mid\mult_PL\ges k\},\\ L^{[(k)]}&:=\{P\in L\mid\mult_PL\in(k)\}\q\q(k\in\Z_{>0}),\endaligned$$
with $(k):=k\1\Z$ (similarly for $L_J^{[k]}:=(L_J)^{[k]}$, etc.)
\sk
Note that conditions~(a), (b) are essentially equivalent to the condition in \cite{ESV}. In the case $m\eq 2$, they are equivalent to each other replacing $J$ with its complement.
\sk
As a corollary of Theorem~1, we get the following.
\msn
{\bf Corollary~1.} {\it In the notation of Theorem~$1$, the eigenspace $H^1(\Ff,\C)_{\la}$ with $\ord\,\la\eq m$ is ESV-calculable if one of the following two conditions is satisfied\,$:$}
$$L^{[\ges 3]}\cap L^{[(m)]}\cap {\rm Sing}\,L_J=\emptyset,
\leqno{\rm (a)'}$$
\vskip-7mm
$$\exists\,P\in L\,\,\,\,\,\h{with}\,\,\,\,\,\tfrac{d}{m}\les\mult_PL\,\notin\,m\1\Z.\leqno{\rm (a)''}$$
\msn
{\bf Remark~1.} Condition~(a)$'$ implies (a), and (a)$''$ implies (a)$'$. (For the latter assertion, we can take $J$ so that all the $L_i$ ($i\in J$) pass through $P$, since $|J|=\tfrac{d}{m}$.)
\msn
{\bf Remark~2.} Condition~(a) is equivalent to (a)$'$ if $L^{[(m)]}=L^{[m]}$. So condition~(a) is a rather strong condition in this case (especially when $m$ is small).
\ms
Condition (a)$''$ in Corollary~1 means that we may have the ESV-calculability easily when $L$ contains a point of very big multiplicity.
\sk
It is not necessarily easy to see whether condition~(a) or (b) is satisfied for each example, since we have to verify it {\it for any $\1J\subset\{1,\dots,d\}$ with\1} $|J|\eq d/m$.
These conditions are often unsatisfied if complexity or efficiency of the combinatorics of arrangement is rather high as in the case of reflection arrangements. In order to measure these properties, we introduce the notions of $m$-complexity and $m$-efficiency for $m\gess 3$. (The case $m\eq 2$ cannot be treated similarly.) The $m$-{\it complexity\1} is defined by
$$C_{L,m}:={\max}_P\,n_P\q\q(m\gess 3).$$
Here $P$ runs over points of $L^{[(m)]}$, and $n_P$ is the number of lines in $L$ containing $P$ and another point of $L^{[(m)]}$. This can be defined also as the maximal number of edges containing one vertex of the associated $m$-graph (which will be defined after Problem~1 below). We have the following (see (1.3) below).
\msn
{\bf Theorem~2.} {\it Condition~{\rm(b)} holds if
$$C_{L,m}\les\lceil m/2\rceil=\bl[(m{+}1)/2\br]\q\q(m\gess 3),
\leqno{\rm(c)}$$
with $\lceil\al\rceil:=\min\{k\in\Z\mid k\ges\al\}$ for $\al\in\R$.}
\ms
If $m\eq 3$, condition~(c) is equivalent to that its $m$-graph is {\it unsaturated\1} (since $\lceil m/2\rceil\eq m{-}1$), see the definition after Problem~1 below.
This bound is sharp when $m\eq 3$. Indeed, there are many ESV-non-calculable examples with $C_{L,3}\eq 3$, see Examples~(3.1) and (3.2) below. It is unclear for $m\ges 4$, see Remark~(1.3)\,(iii) for $m\eq4$.
\sk
The $m$-{\it efficiency\1} is the sum of the numbers of {\it local\1} irreducible components of $L$ passing through the points of $L^{[(m)]}$, which is divided by $d$, that is,
$$E_{L,m}:=\msum_{k\ges 1}\,|L^{[km]}|\1\tfrac{km}{d}\q\q(m\ges 3).$$
This measures how many points of $L^{[(m)]}$ lie on one line {\it in average\1} (since each point of $L^{[km]}$ {\it contributes to\1} $k\1m$ lines). This number can become rather large as $\tfrac{d}{m}$ increases, see Examples~(3.3) and (3.6) below.
The number is closely related to condition~(b) by the following.
\msn
{\bf Remark~3.} Condition~(b) is trivially satisfied if $E_{L,m}\less 1$, that is, if $\,\msum_k\,|L^{[km]}|\1k\les\tfrac{d}{m}\eq|J|$.
\ms
It is rather surprising that we cannot improve this easy bound without imposing any additional hypotheses as in Problem~1 below. In most cases, condition~(b) is satisfied if $E_{L,m}\less 2$, although there are certain exceptional cases, especially when $m\eq 3$.
\msn
{\bf Problem~1.} Is condition (b) always satisfied when $E_{L,m}\less2$ with $m\gess4$, assuming that the $m$-graph of $L$ is connected and unsaturated?
\ms
Here the $m$-{\it graph\1} of $L$ can be described as follows: Its vertices are identified with the points of $L^{[(m)]}$, and there is an edge between two vertices if and only if there is a line in $L$ containing them, see (2.1) below for a more precise definition. It is a {\it weighted\1} graph, where a point has weight $k$ if it belongs to $L^{[km]}$. The $m$-graph is called {\it reduced\1} if the weight is 1 for any vertex, or equivalently, if $L^{[(m)]}=L^{[m]}$. In this case, $L$ is called $m$-{\it reduced.} The $m$-graph is called \h{\it unsaturated\1} if the number of edges containing each vertex is at most $m{-}1$, or equivalently, if $C_{L,m}\less m{-}1$. Here we do {\it not\1} assume that the $m$-graph is {\it reduced.}
\msn
{\bf Remark~4.} In the case $m\eq 3$, we have a positive answer to Problem~1 assuming only its last hypothesis (without assuming $E_{L,m}\less 2$, etc.), since the latter condition is equivalent to the hypothesis of Theorem~2 (that is, $C_{L,m}\les\lceil m/2\rceil$). Note that there are examples such that $E_{L,3}>2$, its $3$-graph is reduced, and the last two hypotheses in Problem~1 are satisfied, see Remark~(1.3)\,(ii) below.
\msn
{\bf Remark~5.} When $E_{L,m}$ is more than 2, it is usually difficult to satisfy condition~(a) or (b). However, the situation does not seem quite simple as is seen below.
\msn
(i) If $f=(x^a{-}y^a)(x^a{-}z^a)(y^a{-}z^a)$ ($a\gess 2$) with $m\eq 3$, we have $\tfrac{d}{3}\eq a$, $\,E_{L,3}\eq a{+}1$ or $a$, and $H^1(\Ff,\C)_{\la}$ for $\ord\,\la\eq 3$ is either ESV-non-calculable or calculable, all depending on whether $a\in3\1\Z$ or not, see Example~(3.1) below.
\msn
(ii) For $m\eq 3$, there are two line arrangements such that $E_{L,3}\eq\tfrac{d}{3}\eq 3$ for both, but $H^1(\Ff,\C)_{\la}$ with $\ord\,\la\eq 3$ is either ESV-calculable or not, see Example~(3.2) below.
\msn
(iii) For $m\eq 3$, there is a family of line arrangements parametrized by $a\in\Z_{>3}$ such that $\,E_{L,3}=(a{+}3)/2\,$ if $\,a\in3\1\Z,\,$ and $\,(a{+}1)/2\,$ otherwise, where $d\eq3a$. It is rather surprising that $H^1(\Ff,\C)_{\la}$ with $\ord\,\la\eq 3$ {\it is\1} ESV-calculable for {\it any\1} $a>3$.
\msn
(iv) In the Hessian arrangement case, we have $E_{L,4}\eq\tfrac{d}{4}\eq 3$, and $H^1(\Ff,\C)_{\la}$ for $\ord\,\la\eq 4$ {\it is\1} ESV-calculable, see Example~(3.4) below.
\msn
(v) For every $m\gess 3$, there is a projective line arrangement such that $E_{L,m}\eq3$, $\tfrac{d}{m}\eq m$, and $H^1(\Ff,\C)_{\la}$ with $\ord\,\la\eq m$ is ESV-calculable if and only if $m$ is odd, see Example~(3.5) below.
\msn
(vi) In the case of a general hyperplane section of the reflection arrangement of type $G_{31}$, we have $E_{L,3}\eq19$, $E_{L,6}\eq3$, $d\eq60$, and $H^1(\Ff,\C)_{\la}$ for $\ord\,\la\eq 3$ or 6 is never ESV-calculable, see \cite{BDY}, \cite{ac} and Example~(3.6) below.
\ms
Note finally that we have a fundamental problem as follows:
\skn
{\bf Problem~2.} Is there an example such that $H^1(\Ff,\C)_{\la}\ne 0$ with $m\eq\ord\,\la\notin\{1,3\}$ except the Hessian arrangement in (iv) above? (Perhaps one could conjecture that there would be no such example; in particular, $\#\bl\{\la\in\C^*\1{\setminus}\1\{1\}\,\big|\,H^1(\Ff,\C)_{\la}\ne 0\br\}\les 3$ for any line arrangement in $\PP^2$, see also \cite[Conjecture 1.9]{PS}, where an analogue of non-existence of 4-nets except the Hessian is not considered so that the above inequality is conjectured with 3 replaced by 5.)
\sk
This problem is closely related to {\it non-existence of $k$-multinets\1} for $k>4$ together with conjectural one for $k\eq 4$ except the Hessian arrangement \cite[Conjecture 4.3]{Yu2} (and also the non-surjectivity of (1.1.6) below). Recall that non-vanishing Milnor fiber cohomology groups with non-unipotent monodromy are closely related to {\it resonance varieties\1} which are the non-vanishing loci of the corresponding Aomoto complex cohomology, and moreover there is a close relation between resonance varieties and multinets, see \cite{FY}, \cite{LY}, \cite{PY}, \cite{Yu1}, \cite{Yu2}.
It is, however, rather unclear whether these can imply that $H^1(\Ff,\C)_{\la}\eq0$ in the ESV-calculable case with $m\eq\ord\,\la>4$ (since the conditions for multinet do not seem quite strong), see (2.5) below.
\sk
In Section~1 we recall some basics of Aomoto complexes, and prove Theorems~1 and 2.
In Section~2 we explain why the assumptions are needed in Problem~1.
In Section~3 we calculate some examples.
\sk
This work was partially supported by JSPS Kakenhi 15K04816.
\bs\bs
\vbox{\centerline{\bf 1. ESV-calculation}
\bsn
In this section we recall some basics of Aomoto complexes, and prove Theorems~1 and 2.}
\msn
{\bf 1.1.~Aomoto complexes.} Let $L\subset Y:=\PP^{n-1}$ be a reduced hyperplane arrangement of degree $d$. Let $\Sc(L)$ be the intersection poset consisting of any intersections of the irreducible components $L_i$ of $L$ ($i\in\{1,\dots,d\})$. This contains the ambient space $Y$, but not the empty set. Set $U:=Y\setminus L$. By \cite{Br}, \cite{OS}, there is an isomorphism of $\C$-algebras
$$A_{\Sc(L)}^{\ssb}\simto H^{\ssb}(U,\C),
\leqno(1.1.1)$$
where $A_{\Sc(L)}^{\ssb}$ is a quotient algebra of the exterior algebra
$\bigwedge^{\ssb}\bigl(\bigoplus_{i=1}^{d-1}\C e_i\bigr)$ divided by an ideal $\I$, and is called the Orlik-Solomon algebra. Note that the induced affine arrangement on $\C^{n-1}=\PP^{n-1}\setminus L_d$ is used here so that the $e_i$ are identified with $\ddd g_i/g_i$, where $g_i$ is a linear function with a constant term defining $L_i\setminus L_d\subset\C^{n-1}$.
Moreover, the ideal $\I$ is determined by the combinatorial data, see {\it loc.\,cit.}
\sk
From now on, {\it assume\1} $n\eq 3$. Let $\al_i\in\C$ $(i\in\{1,\dots,d\})$ satisfying the conditions\,:
$$\aligned&\q\q\al_i\notin\Z_{>0},\q\q\q\msum_{i=1}^d\,\al_i=0,\\&\al_P:=\msum_{L_i\ni P}\,\al_i\,\notin\,\Z_{>0}\q\q(\forall\,P\in L^{[\ges 3]}),\endaligned
\leqno(1.1.2)$$
where $L^{[\ges 3]}$ is as in the introduction. Set
$$\omega=\msum_{i=1}^{d-1}\al_ie_i\in\A^1_{\Sc(L)}.
\leqno(1.1.3)$$
This defines a complex $(\A^{\ssb}_{\Sc(L)},\omega\wedge)$,
called the {\it Aomoto complex\1} associated to $\omega$.
\sk
Since the $e_i$ are identified with $\ddd g_i/g_i$, we get also a regular singular connection $\nabla^{\omega}$
on $\OO_U$ such that
$$\nabla^{\om}h=\ddd h+h\1\om\q(h\in\OO_U).
\leqno(1.1.4)$$
The main theorem of \cite{ESV} asserts that if condition~(1.1.2) is satisfied, then we have the isomorphisms
$$H^j(\A^{\ssb}_{\Sc(L)},\omega\wedge)\simto H^j_{\DR}\bl(U,(\OO_U,\nabla^{\omega})\br)\q(j\in\Z).
\leqno(1.1.5)$$
In this case we say that the de Rham cohomology of $\nabla^{\om}$ is ESV-calculable.
If the last condition of (1.1.2) is not satisfied, we have only the injectivity of
$$H^1(\A^{\ssb}_{\Sc(L)},\omega\wedge)\into H^1_{\DR}\bl(U,(\OO_U,\nabla^{\omega})\br).
\leqno(1.1.6)$$
In the Milnor fiber cohomology case, no example of a non-surjective morphism seems to be known for $m\,{>}\,3$, see Example~(3.1) below for the case $m\eq 3$.
\msn
{\bf 1.2.~Proof of Theorem~1.} It is well-known (see \cite{Di1}, \cite{CS}, \cite{BS}, \h{etc.}) that the monodromy eigenspace $H^j(\Ff,\C)_{\la}$ vanishes unless $\la^d\eq 1$, and there are local systems $\Lc^{(k)}$ of rank 1 on $U$ ($k=0,...,d{-}1$) such that
$$H^j(\Ff,\C)_{\la}=H^j(U,\Lc^{(k)})\q\q\bl(\la=\exp(-2\pi\sqrt{-1}\1k/d)\br).
\leqno(1.2.1)$$
Moreover the monodromy around any irreducible component of $L$ is given by multiplication by $\la^{-1}=\exp(2\pi\sqrt{-1}\1k/d)$. In particular, $\Lc^{(0)}=\C_U$ so that $H^j(\Ff,\C)_1=H^j(U,\C)$.
\sk
It is also known that the dimensions of $H^j(\Ff,\overline{\Q})_{\la}$ (and $H^j(\Ff,\C)_{\la}$) are {\it independent\1} of $\la\in\mu_d$ with $\ord\,\la\eq m$ for a fixed $m$. (Here $\mu_d:=\{\la\in\C\mid\la^d\eq 1\}$.) Indeed, these $\la$ are called {\it primitive roots of unity\1} of order $m$, and are {\it conjugate\1} to each other under the action of the Galois group of $\overline{\Q}/\Q$ (by the irreducibility of cyclotomic polynomials).
In particular, we may assume $\la\eq\exp(\pm 2\pi\sqrt{-1}/m)$ in (1.2.1), that is, $$\tfrac{k}{d}=\tfrac{1}{m}\q\h{or}\q 1\mi\tfrac{1}{m}.
\leqno(1.2.2)$$
Here one can take a convenient choice as one likes.
\sk
The local system $\Lc^{(k)}$ is then isomorphic to the solution local system of the connection $\nabla^{\om}$ on $\OO_U$ in (1.1) by setting
$$\al_i=\begin{cases}1\mi\tfrac{1}{m}&(i\in J),\\-\tfrac{1}{m}&(i\notin J),\end{cases}\q\q\h{or}\q\q\al_i=\begin{cases}-1\pl\tfrac{1}{m}&(i\in J).\\ \tfrac{1}{m}&(i\notin J),\end{cases}
\leqno(1.2.3)$$
where $J\subset I:=\{1,\dots,d\}$ is a subset with $|J|=\tfrac{d}{m}$. (Note that a local system of rank 1 on the complement of a curve in $\PP^2$ is determined by the local monodromies around the {\it smooth points\1} of the curve, taking the tensor product of a local system with the inverse of another, since the singular points of the curve have codimension 2 in $\PP^2$, and $\PP^2$ contains a simply connected Zariski-open subset $\C^2$.)
\sk
In order to get the isomorphism (1.1.5) (that is, the $\la$-eigenspace of the Milnor fiber cohomology is ESV-calculable), the last condition of (1.1.2) should be satisfied for the $\al_i$ defined by one of the equalities in (1.2.3). Setting
$$J_P:=\bl\{i\in J\,\big|\,L_i\ni P\br\},\q I_P:=\bl\{i\in I\,\big|\,L_i\ni P\br\},$$
the conditions are expressed respectively by the following inequalities for $P\in L^{[(m)]}\cap L^{[\ges 3]}$\,:
$$m\1|J_P|\les|I_P|,\q m\1|J_P|\ges|I_P|.$$
Let $k_P,k'_P\in\Z_{>0}$ with $P\in L_J^{[k_P]}\cap L^{[k'_Pm]}$, that is, $|J_P|=k_P$, $|I_P|=k'_Pm$. The above conditions are then equivalent to conditions~(a) and (b) respectively. This finishes the proof of Theorem~1.
\msn
{\bf Remark~1.2.} It may be possible to calculate $H^1(\Ff,\C)_{\la}$ using another method as follows: Let $\pi:\Yt\to Y:=\PP^2$ be a blow-up along a sufficiently general point $P\in L_d$ with $\jt:U\into\Yt$ a natural inclusion. Set $\Ft:=\RR\jt_*\Lc^{(k')}[2]$ for $k'\eq d/m\,$ or $\,d\mi d/m$. This is a direct factor of the underlying $\C$-complex of a mixed Hodge module, and has the weight filtration $W$ such that $\Gr^W_2\Ft={\rm IC}_{\Yt}\Lc^{(k')}$ and $\Gr^W_i\Ft=0$ for $i\ne 2,3$. (In the case $\pi(y)\in L^{(km)}$ ($k\in\Z_{>0}$), we have $\,\dim_{\C}\Hc^j\Gr^W_i\!\Ft_y=km\mi 2\,$ if $\,j\eq i\mi1\eq 1$ or $2$, and it is zero otherwise.)
\sk
We have a natural projection $\rho:\Yt\to C:=\PP^1$ so that $\Yt$ is a $\PP^1$-bundle over $C$, and the exceptional divisor of the blow-up is the zero-section (where $C$ is identified with the set of lines on $\PP^2$ passing through $P$). In the notation of \cite{BBD}, we have
$$^{\bf p}\Hc^j\RR\rho_*\Gr^W_i\!\Ft=0\q\q(j\ne0,\,\,i\eq2,3).$$
Set $\G:=\RR\rho_*\Ft\,(={}^{\bf p}\Hc^0\RR\rho_*\Ft)$. This has the weight filtration $W$ such that
$$\Gr^W_i\G={}^{\bf p}\Hc^0\RR\rho_*\Gr^W_i\!\Ft\q\q(i\eq2,3).$$
Set $C':=C\setminus\rho\bl(\pi^{-1}({\rm Sing}\,L)\br)$. We see that ${\rm Supp}\,\Gr^W_3\G\subset C\setminus C'$, and $\Gr^W_2\G$ is an intersection complex whose restriction to $C'$ is a shifted local system of rank $d{-}2$. (The latter can be calculated by using the direct image from the blowing-up of $\Yt$ at the non-normal crossing singular points of $L$.) We then get that
$$H^1(\Ff,\C)_{\la}\cong H^{-1}(C,\G)=\Gamma(C',\G[-1]|_{C'}),$$
where $\la=\exp(\pm 2\pi\sqrt{-1}/m)$. (Note that the monodromy is defined over $\Q$.) It might be possible to obtain the last term of the equalities using a Picard-Lefschetz type formula and calculating the vanishing cycles of $\G$ at each point of $C\setminus C'$.
\msn
{\bf 1.3.~Proof of Theorem~2.} For $P\in L^{[(m)]}$, let $r_P$ be the number of lines $L_i\subset L$ such that
$$L_i\cap L^{[(m)]}\eq\{P\}.$$
Let $k_P\in\Z_{>0}$ such that $P\in L^{[k_Pm]}$. Condition~(c) means that
$$n_P=k_Pm\mi r_P\les\lceil m/2\rceil,\q\h{that is,}\q r_P\ges k_Pm\mi(m{+}1)/2,
\leqno(1.3.1)$$
since $\lceil m/2\rceil=\bl[(m{+}1)/2\br]$. By Remark~(1.3)\,(i) below, the assertion can be reduced to the $m$-reduced case, and we may assume $L$ is $m$-reduced.
\sk
In this case, it is enough to show the following.
\msn
\rlap{\rm(1.3.2)}\hskip1.7cm\hangindent=1.7cm
For any $m$-reduced line arrangement $L$ with condition~(c) satisfied, there are $r$ lines $L_i\subset L$ ($i\in\{1,\dots,r\}$) with $m\1r\less d$ and $L^{[m]}\subset\mcup_{i=1}^r\,L_i$, where $d:=\deg L$.
\msn
Here we do not assume $d/m\in\Z$.
\sk
We proceed by induction on $|L^{[m]}|$. The assertion holds in the case $|L^{[m]}|\eq 1$, since $d\gess m$. Assume $|L^{[m]}|\gess 2$. Take any line $L_1\subset L$ such that $|L_1\cap L^{[m]}|\gess 2$. (If there is no such line, the assertion is reduced to the case $|L^{[m]}|\eq 1$.) Let $L''\subset L$ be the union of $L_1$ and the lines $L_i\subset L$ such that
$$L_1\supset L_i\cap L^{[m]}\ne\emptyset.
\leqno(1.3.3)$$
Let $L'\subset L$ be the union of the other lines so that $L\eq L'\cup L''$ and $\dim L'\cap L''\eq 0$. By condition~(1.3.1) (with $k_P\eq 1$), we have
$$d'':=\deg L''\ges 1+2\bl(m\mi(m{+}1)/2\br)=m.
\leqno(1.3.4)$$
Hence $d''\gess m\1r''$ with $r'':=1$. By inductive hypothesis, there are $r'$ lines in $L'$ covering $L'{}^{[m]}$ with $d':=\deg L'\gess m\1r'$. By (1.3.3) we have $L'{}^{[m]}\eq L^{[m]}\setminus L_1$. Hence there are $r$ lines in $L$ covering $L^{[m]}$ with $r:=r'\pl r''$.
Theorem~2 then follows, since $d\eq d'\pl d''$.
\msn
{\bf Remark~1.3}\,(i). For a line arrangement $L$ with $C_{L,m}\less m$, we have the decomposition $L\eq L'\cup L''$ such that $L'$ is $m$-reduced (that is, $L'{}^{[(m)]}\eq L'{}^{[m]}$), $L''{}^{[km]}\eq L^{[(k{+}1)m]}$ ($\forall\,k\in\Z_{>0}$), and $|L''_i\cap L''{}^{(m)}|\eq1$ for any line $L''_i\subset L''$. This $L'$ is called the $m$-reduced arrangement associated to $L$. It is easy to see that condition~(a) (resp.~(b)) is satisfied for $L$ if and only if it is satisfied for $L'$.
\msn
{\bf Remark~1.3}\,(ii). For $m\eq 3$, there are many examples such that the $3$-graph is connected, reduced, and unsaturated with $E_{L,3}>2$. For instance, take sufficiently general $n$ lines $L_i\in\PP^2$ ($i\in\{1,\dots,n\}$) such that their union $L':=\mcup_{i=1}^n\,L_i$ is a divisor with normal crossings. For each singular point of $L'$, choose a sufficiently general line passing through it so that its intersection with ${\rm Sing}\,L'$ consists only of this singular point and we have $L^{[(3)]}=L^{[3]}={\rm Sing}\,L'$, where $L$ is the union of these lines and $L'$. Then
$$\big|L^{[(3)]}\big|=\big|L^{[3]}\big|=n(n{-}1)/2,\q d:=\deg L=n(n{+}1)/2.$$
Hence
$$E_{L,3}=3(n{-}1)/(n{+}1)>2\q\h{if}\q n\gess 6.$$
\msn
{\bf Remark~1.3}\,(iii). Put $m':=m{-}1$. For a line arrangement $L'$ with $L'{}^{[(m)]}\eq \emptyset$, there is a line arrangement $L$ whose $m$-graph is the same as the $m'$-graph of $L'$ and
$$d\eq d'\pl\msum_{k\ges1}\,|L'{}^{[km']}|\1k,$$
with $d:=\deg L$, $d':=\deg L'$, by adding sufficiently general $k$ lines passing through each point of $L'{}^{[km']}$ for $k\ges 1$. (Here we do {\it not\1} assume that $L$ is $m$-reduced.) We have
$$d/m\eq\bl(d'\pl\msum_k\,|L'{}^{[m']}|\1k\br)/m=(d'/m')\bl(m'\pl E_{L',m'}\br)/m.$$
This means that $d/m$ is close to $d'/m'$ if $E_{L',m'}$ is close to 1. It does not seem, however, that this argument can be used to show that Theorem~2 is sharp for $m\eq 4$. Here we can apply the construction in (2.3) below for $L'$, and consider the union of $\gamma_j(L)$ $(j\in\{1,\dots,m\}$), where the $\gamma_j\in{\rm Aut}(\PP^2)$ fix some appropriate smooth point of $L$, and are sufficiently general among the automorphisms satisfying this condition. However, in the case $L'$ is obtained by applying (2.3) below to the second arrangement in Example~(3.2) below with $m'\eq 3$, for instance, we get that $d'/m'\eq 3\pl k_{L''}$, $d/m\eq 19/4\pl k_{L''}$ with $k_{L''}$ as in (2.3) below.
\bs\bs
\vbox{\centerline{\bf 2. Necessity of the hypotheses in Problem~1.}
\bsn
In this section we explain why the assumptions are needed in Problem~1.}
\msn
{\bf 2.1.~$m$-graphs.} For a line arrangement $L\subset\PP^2$ and $m\in\Z_{\ges3}$, the $m$-graph is defined as follows: Its vertices are identified with the points of $L^{[(m)]}$. Its edges correspond to lines in $L$ containing at least two points of $L^{[(m)]}$, and are expressed by closed connected smooth real curves whose both ends are vertices. An edge may contain vertices in its interior so that the vertices contained in it are identified with the points of $L^{[(m)]}$ contained in the corresponding line in $L$. We assume that two edges containing a same vertex have different limit tangent lines. (Note that an intersection of edges is not always a vertex of the graph unless a vertex is marked at the intersection point. If two edges intersect at a point which is not a vertex, we consider that they do not intersect as if the graph is in $\R^3$.)
\sk
This is a weighted graph, where a vertex has weight $k$ if its corresponding point belongs to $L^{[km]}$. We can decide whether condition~(a) or (b) is satisfied by looking at this graph. Indeed, the number of lines in $L$ containing only a given vertex is the difference between the weight of the vertex multiplied by $m$ and the number of edges containing the vertex. The weights are omitted if these are always 1, that is, if $L$ is $m$-reduced.
(Note that the $m$-graph cannot be used to determine whether $L$ supports a multinet. For this we have to employ the {\it complete graph\1} whose vertices are all the singular points of $L$.)
\sk
For instance, let $f=(x^a\mi y^a)(x^a\mi z^a)(y^a\mi z^a)$ ($a\eq 2,3$) as in Example~(3.1) below. For $a\eq 2$ (resp.\ 3), this is the simplest example with $H^1(\Ff,\C)_{\la}\ne 0$ for $\la\ne 1$ (resp.\ one of the simplest ESV-non-calculable examples).
Its $3$-graph can be drawn as follows (here we use piecewise smooth curves whose tangent lines vary continuously by technical reasons):
\par\q\q\q\q\q\q\q
\raise7mm\h{$\setlength{\unitlength}{7mm}
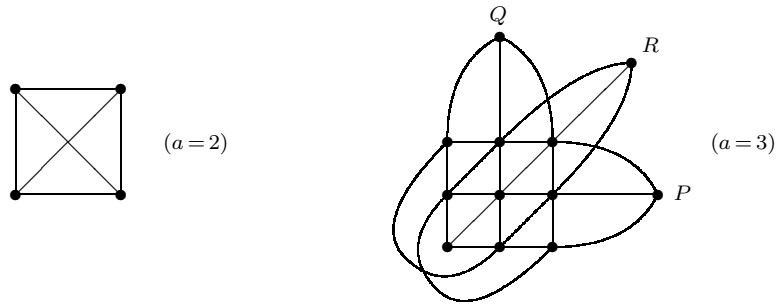
\begin{picture}(4,5)
\multiput(1,1)(2,0){2}{\line(0,1){2}}
\multiput(1,1)(0,2){2}{\line(1,0){2}}
\multiput(1,1)(2,0){2}{\circle*{.2}}
\multiput(1,3)(2,0){2}{\circle*{.2}}
\put(1,1){\line(1,1){2}}
\put(1,3){\line(1,-1){2}}
\end{picture}$}\!\!\!\!
\raise20mm\h{$\scriptstyle(a\eq 2)$}
\q\q\q\q\q
\h{$\setlength{\unitlength}{7mm}
\begin{picture}(6,6)
\multiput(1,1)(1,0){3}{\line(0,1){2}}
\multiput(1,1)(0,1){3}{\line(1,0){2}}
\multiput(1,1)(1,0){3}{\circle*{.2}}
\multiput(1,2)(1,0){3}{\circle*{.2}}
\multiput(1,3)(1,0){3}{\circle*{.2}}
\put(5,2){\circle*{.2}}
\put(2,5){\circle*{.2}}
\put(4.5,4.5){\circle*{.2}}
\put(5.3,1.9){$\scriptstyle P$}
\put(1.8,5.3){$\scriptstyle Q$}
\put(4.7,4.7){$\scriptstyle R$}
\qbezier(1,3)(1,4.5)(2,5)
\qbezier(2,3)(2,4)(2,5)
\qbezier(3,3)(3,4.5)(2,5)
\qbezier(3,1)(4.5,1)(5,2)
\qbezier(3,2)(4,2)(5,2)
\qbezier(3,3)(4.5,3)(5,2)
\qbezier(3,2)(4.5,3.5)(4.5,4.5)
\qbezier(2,1)(2.5,1.5)(3,2)
\qbezier(.2,.8)(1,0)(2,1)
\qbezier(.2,.8)(-.5,1.5)(1,3)
\qbezier(2,3)(3.5,4.5)(4.5,4.5)
\qbezier(1,2)(1.5,2.5)(2,3)
\qbezier(.8,.2)(0,1)(1,2)
\qbezier(.8,.2)(1.5,-.5)(3,1)
\put(1,1){\line(1,1){3.5}}
\end{picture}$}\!\!
\raise20mm\h{$\scriptstyle(a\eq 3)$}
\skn
These pictures are based on the restriction of $L$ to the affine space $\C^2\subset\PP^2$ defined by the equation $(x^a\mi y^a)(x^a\mi 1)(y^a\mi 1)=0$. For $a\eq 3$, the nine vertices written in the left lower part correspond to
$$(\zeta^i,\zeta^j)\in\C^2\,\,\,\,\bl((i,j)\in(\Z/3\1\Z)^2\br),$$
with $\zeta=e^{2\pi\sqrt{-1}/3}$. Among the remaining three vertices, $P,Q$ correspond to the points at infinity, and $R$ corresponds to the origin of $\C^2$.
(One can find similar pictures in the literature. Note that ${\rm Sing}\,L=L^{[3]}$ if $a\eq 3$.) In the case of {\it real\1} line arrangements, one can draw its $m$-graph by modifying its restriction to $\R^2$.
\msn
{\bf 2.2.~Disconnected $m$-graph case.} Assume $L=L'\cup L''$ with $d'$, $d''$ divisible by $m$, $\,d\eq d'\pl d''$ (where $d':=\deg L'$, etc.), and moreover
$$L^{[(m)]}=L'{}^{[(m)]}\sqcup L''{}^{[(m)]}.$$
\vskip-3mm\nin
Then
\vskip-5mm
$$E_{L,m}\eq\tfrac{d'}{d}\1 E_{L',m}+\tfrac{d''}{d}\1 E_{L'',m},$$
(in particular, $E_{L',m}>E_{L,m}>E_{L'',m}$ if $\,E_{L',m}>E_{L'',m}$).
This implies a negative answer to Problem~1 in the case the $m$-graph is disconnected. (Note that $E_{L'',m}\eq 1$ if $L''$ is non-essential, that is, if $d''\eq km$ and $L''{}^{[km]}\ne\emptyset$ with $k\in\Z_{>0}$.)
\msn
{\bf 2.3.~Subarrangements of $m$-star type, I.} For any line arrangement $L'$, any $k_0\in\Z_{>0}$, and any line $L_0\subset L'$ containing $P\in L'{}^{[k_Pm]}$ with $k_P\in\Z_{>0}$, we can construct a line arrangement $L=L'\cup L''$ such that $L''$ has $m$-graph of $m$-star type with following conditions satisfied\1:
\skn
\rlap{\rm(i)}\hskip.8cm\hangindent=.8cm
The intersection $L'\cap L''$ coincides with the union of lines in $L'$ containing $P$.
\skn
\rlap{\rm(ii)}\hskip.8cm\hangindent=.8cm
$\,L'{}^{[(m)]}\cup L''{}^{[(m)]}=L^{[(m)]},\q L'{}^{[(m)]}\cap L''{}^{[(m)]}=\{P\},$
\skn
\rlap{\rm(iii)}\hskip.8cm\hangindent=.8cm
There is $Q\in L_0\cap L''{}^{[k_0m]}$ such that $|L''_i\cap L''{}^{[(m)]}|=2$ for any line $L''_i\subset L''$ with $Q\in L''_i$.
\skn
\rlap{\rm(iv)}\hskip.8cm\hangindent=.8cm
$\,\deg L''=\msum_{k\ges 1}\,|L''{}^{[km]}|\1km-k_Qm,\q |L''{}^{[(m)]}|\eq k_Qm{+}1\q(k_Q:=k_0).$
\msn
(We say that $L''$ has $m$-graph of $m$-star type if conditions (iii) and (iv) are satisfied with $L_0$ forgotten.)
\sk
Put $d:=\deg L$, $\,d':=\deg L'$, $\,d'':=\deg L''$. We have
$$\aligned&d\mi d'=d''\mi k_Pm=k_{L''}\1m,\q\q E_{L'',m}=(k_{L''}{+}k_P{+}k_Q)\1m/d'',\\&\h{with}\q\q\q k_{L''}:=\msum_{k\ges 1}\,|L''{}^{[km]}|\1k-(k_P\pl k_Q).\endaligned$$
Setting
\vskip-5mm
$$\Et_{L'',m}:=(k_{L''}\pl k_Q)\1m/(d\mi d')=1+k_Q/k_{L''},$$
we then get that
$$E_{L,m}\eq\tfrac{d'}{d}\1 E_{L',m}+\tfrac{d-d'}{d}\,\Et_{L'',m}.$$
\skn
Here we can make $E_{L,m}$ arbitrarily close to 1, since $k_{L''}$ can be arbitrarily large. (We do not assume that $L''$ is $m$-reduced.) This implies a negative answer to Problem~1 unless we assume the $m$-graph is unsaturated.
Note that if there is $J\subset\{1,\dots,d\}$ such that condition~(a) (resp.\ (b)) is satisfied, then $|J\cap\{d'{+}1,\dots,d\}|$ must be at most (resp.\ at least) $(d{-}d')/m=k_{L''}$.
(Note also that we can arrange so that $C_{L'',m}\eq m$, if $k_0\eq 1$. Indeed, $L''{}^{[k'']}$ can be non-empty for some very large $k''$, since we do not assume that $L''$ is $m$-reduced. This is needed for Remark~(1.3)\,(iii).)
\msn
{\bf 2.4.~Subarrangements of $m$-star type, II.} In the case $m\eq 3$, it is not difficult to construct an arrangement having a subarrangement with $m$-graph of $m$-star type and such that condition~(b) is not satisfied although the $m$-efficiency is at most 2. For instance, consider a line arrangement $L'$ of degree 8 whose $3$-graph is as follows:
$$\setlength{\unitlength}{5mm}
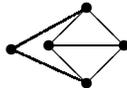
\begin{picture}(5,2.2)
\qbezier(0,.9)(0,.9)(2,2)
\qbezier(0,.9)(0,.9)(2,0)
\put(1,1){\line(1,0){2}}
\put(1,1){\line(1,-1){1}}
\put(1,1){\line(1,1){1}}
\put(3,1){\line(-1,1){1}}
\put(3,1){\line(-1,-1){1}}
\put(0,.9){\circle*{.3}}
\put(1,1){\circle*{.3}}
\put(2,2){\circle*{.3}}
\put(2,0){\circle*{.3}}
\put(3,1){\circle*{.3}}
\end{picture}$$
The number of edges is 7, and there is a line $L'_0$ in $L'$ which does not contribute to an edge. This line passes through the leftest vertex.
\sk
Let $L=\mcup_{i=1}^3\,L^{(i)}$, where $L^{(i)}$ ($i=1,2,3$) is as above (with $\deg L^{(i)}~8$) and $L^{(i)}\cap L^{(j)}$ is discrete ($i\ne j$). Assume
$$\mcap_{i=1}^3\,L_0^{(i)}=\{P\}\not\subset\mcup_{i=1}^3\,{\rm Sing}\,L^{(i)},\q L^{[3]}=\h{$\bigsqcup$}_{i=1}^3\,L^{(i)[3]}\sqcup\{P\},$$
where $L_0^{(i)}\subset L^{(i)}$ is the unique line containing only one point of multiplicity 3 (and not contributing to an edge). Then $\deg L=24$, $E_{L,3}\eq2$, and conditions~(a), (b) are both unsatisfied.
\sk
It is also possible to consider the case $L=\mcup_{i=1}^6\,L^{(i)}$ where $L^{(i)}$ ($i\in\{1,\dots,4\}$) is as above (with $|L^{(i)}|=8$), $|L^{(5)}|=6$, $|L^{(5)[3]}|=4$, $|L^{(6)}|=1$, $L^{(i)}\cap L^{(j)}$ is discrete ($i\ne j$), and
$$\mcap_{i=1}^4\,L_0^{(i)}\cap L^{(5)}\cap L^{(6)}=\{P\}\not\subset\mcup_{i=1}^5{\rm Sing}\,L^{(i)},\q L^{[(3)]}=\h{$\bigsqcup$}_{i=1}^5\,L^{(i)[3]}\sqcup\{P\}.$$
Here $L^{(5)}$ is as in Example~(3.1) below with $a\eq 2$.
Then $\deg L=39$, $|L^{[3]}|=24$, $|L^{[6]}|=1$, and $E_{L,3}=2$. Conditions~(a) and (b) are both unsatisfied.
\sk
One can replace $L^{(5)}$ in the above example by a non-essential arrangement of degree $3a$ ($a\in\Z_{>2}$) so that $|L^{(5)}|\eq 3a$, $L^{(5)}{}^{[3a]}\ne\emptyset$. In this case we have
$$\aligned\deg L&=3(a{+}11),\q|L^{[3]}|=20,\q|L^{[6]}|=1,\q|L^{[3a]}|=1,\\ E_{L,3}&=(a{+}22)/(a{+}11).\endaligned$$
We can verify that conditions~(a) and (b) are both unsatisfied.
\msn
{\bf 2.5.~Relation with non-existence of $k$-multinets for $k>4$.}
It is known that there is a close relation between resonance varieties, pencils, and multinets, and there are no $k$-multinets for $k>4$, see \cite{FY}, \cite{LY}, \cite{PY}, \cite{Yu1}, \cite{Yu2}.
However, it seems unclear whether these imply the vanishing of $H^1(\Ff,\C)_{\la}$ with $m:=\ord\,\la>4$ in the ESV-calculable case.
\sk
Indeed, if this eigenspace does not vanish, then $L$ supports a $k$-multinet, hence there is a partition
$$I:=\{1,\dots,d\}=\h{$\bigsqcup$}_{j=1}^k\,I_j$$
together with multiplicities $m_i\in\Z_{>0}$ ($i\in I$) such that ${\rm GCD}(m_i)=1$, and in the notation of (1.1) we would have the following (see for instance \cite[Theorem 3.7]{Yu2})\,:
$$\om=\msum_{j=1}^k\,c_j\eta^{(j)}\q\q\h{with}\q\q\eta^{(j)}=\msum_{i\in I_j\setminus\{d\}}m_ie_i.
\leqno(2.5.1)$$
Here $c_j\in\C$ ($j\in\{1,\dots,k\}$) with $\msum_{j=1}^k\,c_j=0$, and $\om$ is defined by (1.1.3) and (1.2.3).
The equality (2.5.1) then implies that $c_j\in\Q$, and $J$ is compatible with the partition of $I$, that is, $J$ is a union of $I_j$ ($j\in K$) for some subset $K\subset\{1,\dots,k\}$ (since $m_i\in\Z_{>0}$).
We would then get that the $m_i$ ($i\in I_j$) depend only on $j$ (denoted by $m_j$) using (1.2.3). However, it does not seem quite clear whether the $m_j$ are independent of $j\in\{1,\dots,k\}$. Consider, for instance, the case
$$k\eq3,\q c_1\eq2,\q c_2\eq c_3\eq{-}1,\q |I_1|\eq b,\q |I_2|\eq|I_3|\eq ab,\q\ m\eq 2a{+}1$$
and $m_i=a$ if $i\in I_1$, and 1 otherwise, where $a,b\in\Z_{\ges 2}$ (although it does not seem very clear whether this can really happen), see \cite[Definition 3.5 and Remark 3.6]{Yu2} for multinets.
\bs\bs
\vbox{\centerline{\bf 3. Examples}
\bsn
In this section we calculate some examples.}
\msn
{\bf Example~3.1.} For any integer $a\gess 2$, let
$$f=(x^a{-}y^a)(x^a{-}z^a)(y^a{-}z^a),$$
with $d\eq3\1a$. This is a reflection arrangement of type $G(a,a,3)$, see \cite[p.~280]{OT}. For $a\eq2$, this is the simplest example with $H^1(\Ff,\C)_{\la}\ne 0$ for $\la\ne 1$. We have $L^{[k]}\eq\emptyset$ ($k\ne 3,a$), and
$$|L^{[3]}|\eq a^2,\,\,\,|L^{[a]}|\eq3\,\,\,\,(a\ne 3),\q\,\,\,|L^{[3]}|\eq 12\,\,\,\,(a=3).$$
Hence
$$\aligned\tfrac{d}{3}\eq a,\q&E_{L,3}=\begin{cases}a{+}1&(a\in 3\1\Z),\\a&(a\notin 3\1\Z),\end{cases}\\ \tfrac{d}{a}\eq 3,\q &E_{L,a}=1\q(a\ne 3).\endaligned$$
If $a\notin3\1\Z$, we see that $H^1(\Ff,\C)_{\la}$ for $\ord\,\la\eq3$ is ESV-calculable, where $J$ can be given by any of the three factors of $f$. In the other case, it is not ESV-calculable. Its dimension is either 2 or 1, depending on whether $a\in3\1\Z$ or not, see \cite{Di3}. Note that it is 1 in the ESV-calculable case, since we have a 3-net, see \cite{BDS}, etc. When $a\in3\1\Z$, we then get an example with (1.1.6) non-surjective. (These can be confirmed by using a small computer program calculating the pole order spectral sequence when $a\less 6$.)
\msn
{\bf Example~3.2.} Let
$$\aligned&f_1=xyz(y{+}2z)(x{-}y)(x{-}y{+}z)(x{+}y{-}z)(x{+}y{+}2z)(x{-}2y{-}z),\\&f_2=xyz(x{+}y)(y{+}z)(x{+}3z)(x{+}2y{+}z)(x{+}2y{+}3z)(2x{+}3y{+}3z).\endaligned$$
The pictures of the restrictions of the arrangements to $\C^2=\PP^2\setminus\{z\eq0\}$ are as below. (Note that the line at infinity belongs to $L$, and there are three points of multiplicity 3 at infinity.)
$$\raise8mm\h{$\setlength{\unitlength}{0.5cm}
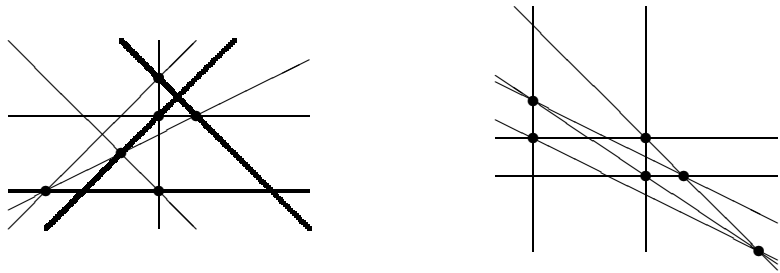
\begin{picture}(8,5)
\linethickness{.1mm}
\put(0,3){\line(1,0){8}}
\put(4,0){\line(0,1){5}}
\put(0,0){\line(1,1){5}}
\put(0,5){\line(1,-1){5}}
\put(0,0.5){\line(2,1){8}}
\linethickness{0.5mm}
\put(0,1){\line(1,0){8}}
\qbezier(1,0)(2,1)(6,5)
\qbezier(3,5)(4,4)(8,0)
\put(4,4){\circle*{.3}}
\put(4,3){\circle*{.3}}
\put(4,1){\circle*{.3}}
\put(5,3){\circle*{.3}}
\put(3,2){\circle*{.3}}
\put(1,1){\circle*{.3}}
\put(4,4){\circle*{.3}}
\end{picture}$}\q\q\q\q\q\q
\h{$\setlength{\unitlength}{5mm}
\begin{picture}(8,8)
\put(4,1){\line(0,1){6.5}}
\put(0,4){\line(1,0){7.5}}
\put(0.5,7.5){\line(1,-1){7}}
\put(0,3){\line(1,0){7.5}}
\put(1,1){\line(0,1){6.5}}
\put(0,5.5){\line(2,-1){7.5}}
\put(0,4.5){\line(2,-1){7.5}}
\put(4,3){\line(3,-2){3.5}}
\put(4,3){\line(-3,2){4}}
\put(7,1){\circle*{.3}}
\put(4,3){\circle*{.3}}
\put(5,3){\circle*{.3}}
\put(4,4){\circle*{.3}}
\put(1,4){\circle*{.3}}
\put(1,5){\circle*{.3}}
\end{picture}$}$$
\vskip-4mm\
Their associated 3-graphs are as below. (To see whether condition~(a) or (b) is satisfied, it seems better to use the associated 3-graphs, since unnecessary information is hidden.)
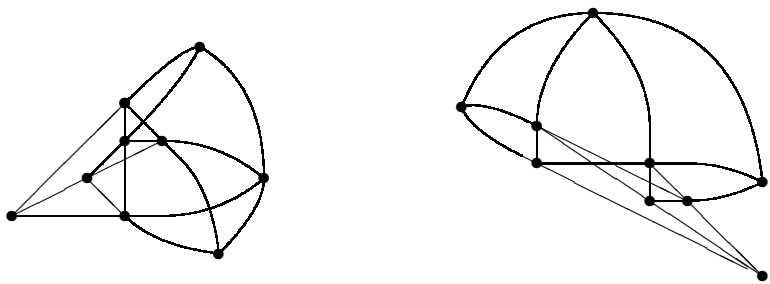
\sk
$$\raise8mm\h{$\setlength{\unitlength}{0.5cm}
\begin{picture}(8,5)
\put(4,3){\line(1,0){1}}
\put(4,1){\line(0,1){3}}
\put(1,1){\line(1,1){3}}
\put(3,2){\line(1,-1){1}}
\put(1,1){\line(2,1){4}}
\put(1,1){\line(1,0){4}}
\qbezier(3,2)(4,3)(4.5,3.5)
\qbezier(4,4)(5,3)(5.5,2.5)
\qbezier(4,4)(5.5,5.5)(6,5.5)
\qbezier(4,3)(5.5,4.5)(6,5.5)
\qbezier(5,1)(6.5,1)(7.7,2)
\qbezier(5,3)(6.5,3)(7.7,2)
\qbezier(4,1)(4.8,.2)(6.5,0)
\qbezier(5.5,2.5)(6.3,1.7)(6.5,0)
\qbezier(7.7,2)(7.7,1.2)(6.5,0)
\qbezier(7.7,2)(7.7,4.5)(6,5.5)
\put(4,4){\circle*{.3}}
\put(4,3){\circle*{.3}}
\put(4,1){\circle*{.3}}
\put(5,3){\circle*{.3}}
\put(3,2){\circle*{.3}}
\put(1,1){\circle*{.3}}
\put(4,4){\circle*{.3}}
\put(7.7,2){\circle*{.3}}
\put(6,5.5){\circle*{.3}}
\put(6.5,0){\circle*{.3}}
\end{picture}$}\q\q\q\q\q\q
\h{$\setlength{\unitlength}{5mm}
\begin{picture}(8,8)
\put(5,3){\line(0,1){2}}
\put(2,4){\line(1,0){4}}
\put(5,4){\line(1,-1){3}}
\put(5,3){\line(1,0){1}}
\put(2,4){\line(0,1){1}}
\put(2,5){\line(2,-1){4}}
\put(1.6,4.2){\line(2,-1){6}}
\put(2,5){\line(3,-2){6}}
\qbezier(2,5)(2,6.5)(3.5,8)
\qbezier(5,5)(5,6.5)(3.5,8)
\qbezier(1.6,4.2)(.2,4.9)(0,5.5)
\qbezier(2,5)(.6,5.7)(0,5.5)
\qbezier(6,3)(7,3)(8,3.5)
\qbezier(6,4)(7,4)(8,3.5)
\qbezier(0,5.5)(1,8)(3.5,8)
\qbezier(8,3.5)(7.5,8)(3.5,8)
\put(8,1){\circle*{.3}}
\put(5,3){\circle*{.3}}
\put(6,3){\circle*{.3}}
\put(5,4){\circle*{.3}}
\put(2,4){\circle*{.3}}
\put(2,5){\circle*{.3}}
\put(3.5,8){\circle*{.3}}
\put(0,5.5){\circle*{.3}}
\put(8,3.5){\circle*{.3}}
\end{picture}$}$$
\vskip-3mm
The first arrangement is a 3-net. This is one of the second simplest examples for non-vanishing Milnor fiber cohomology with non-unipotent monodromy.
The second one is one of the simplest ESV-non-calculable examples, see also \cite{CS}.
In both cases we have $d=9$, and $|L^{[(3)]}|=|L^{[3]}|=9$. Hence
$$E_{L,3}=\tfrac{d}{3}=3.$$
One can verify that $H^1(\Ff,\C)_{\la}$ with $\ord\,\la\eq3$ is ESV-calculable in the first case (where the lines in $J$ are indicated by thick lines), but not in the second case. (Here the argument is divided into two cases depending on whether $J$ contains the divisor at infinity or not.)
\msn
{\bf Example~3.3.} For $a\ges4$, let
$$f=\mprod_{i=0}^{a-1}\1(x\mi iz)\,\mprod_{j=0}^{a-1}\1(y\mi jz)\,\mprod_{k=0}^{a-1}\1(x\pl y\mi kz),$$
with $d\eq3a$. We have $|L^{[3]}|=\tfrac{a(a+1)}{2}$, $|L^{[a]}|=3$, hence
$$E_{L,3}=\begin{cases}(a{+}3)/2&(a\in3\1\Z),\\(a{+}1)/2&(a\notin3\1\Z).\end{cases}$$
It is rather surprising that $H^1(\Ff,\C)_{\la}$ with $\ord\,\la\eq3$ is always ESV-calculable independently of whether $a\in 3\1\Z$ or not. Indeed, if $a\notin3\1\Z$, we can take $J$ corresponding to any one of the three factors of the above factorization. (This is trivial.) In the case $a=3a'$, we can take $J$ corresponding to the polynomial
$$g=\mprod_{i=0}^{a'-1}\1(x\mi iz)\,\mprod_{j=0}^{a'-1}\1(y\mi jz)\,\mprod_{k=2a'}^{a-1}\1(x\pl y\mi kz).$$
The picture in the case $a\eq 6$ is as below.
$$\setlength{\unitlength}{5mm}
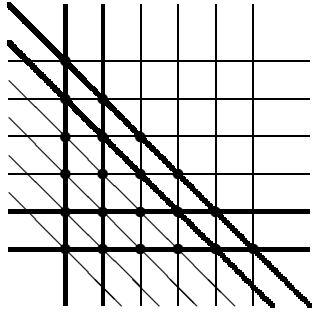
\begin{picture}(8,8)
\linethickness{.1mm}
\multiput(3.5,0)(1,0){4}{\line(0,1){8}}
\multiput(0,3.5)(0,1){4}{\line(1,0){8}}
\put(0,6){\line(1,-1){6}}
\put(0,5){\line(1,-1){5}}
\put(0,4){\line(1,-1){4}}
\put(0,3){\line(1,-1){3}}
\linethickness{.5mm}
\qbezier(0,8)(0,8)(8,0)
\qbezier(0,7)(0,7)(7,0)
\multiput(1.5,0)(1,0){2}{\line(0,1){8}}
\multiput(0,1.5)(0,1){2}{\line(1,0){8}}
\multiput(1.5,6.5)(1,-1){6}{\circle*{.3}}
\multiput(1.5,5.5)(1,-1){5}{\circle*{.3}}
\multiput(1.5,4.5)(1,-1){4}{\circle*{.3}}
\multiput(1.5,3.5)(1,-1){3}{\circle*{.3}}
\multiput(1.5,2.5)(1,-1){2}{\circle*{.3}}
\multiput(1.5,1.5)(1,-1){1}{\circle*{.3}}
\end{picture}$$
Notice that we lose the ESV-calculability in the case $a\in3\1\Z$, if the last product in the definition of $f$ is taken over $k\in\{1,\dots,a\}$ instead of $k\in\{0,\dots,a{-}1\}$.
\msn
{\bf Example~3.4.} Let
$$f=xyz\,\mprod_{i,j=0}^2\,(\zeta^ix\pl\zeta^jy\pl z)=xyz\,\mprod_{k=0}^2\,(x^3\pl y^3\pl z^3\mi 3\1\zeta^kxyz),$$
with $\zeta=e^{2\pi\sqrt{-1}/3}$. This is the {\it Hessian\1} arrangement, which is a unique example of 4-net, see \cite[Ex.~6.30 and p.\,232]{OT}, \cite[Ex.~3.5]{FY}, \cite[Remark 3.3\,(iii)]{BDS}, \cite[Theorem 8.19]{Di2}. It is conjectured that there are no other 4-nets, see \cite{Yu2}. (This arrangement seems to be confused with a different one in some papers.)
Its 4-graph may be drawn as follows.
$$\setlength{\unitlength}{5mm}
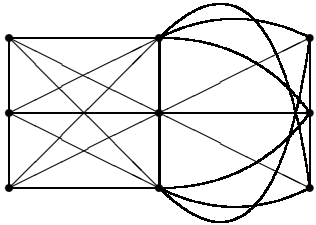
\begin{picture}(8,6.1)
\multiput(0,1)(4,0){3}{\line(0,1){4}}
\multiput(0,1)(0,2){3}{\line(1,0){4}}
\multiput(0,1)(4,0){3}{\circle*{.2}}
\multiput(0,3)(4,0){3}{\circle*{.2}}
\multiput(0,5)(4,0){3}{\circle*{.2}}
\put(0,3){\line(2,-1){4}}
\qbezier(4,1)(6,0)(8,1)
\put(0,5){\line(1,-1){4}}
\qbezier(4,1)(7,-2)(8,5)
\qbezier(4,1)(6.5,1)(8,3)
\put(0,3){\line(2,1){4}}
\qbezier(4,5)(6,6)(8,5)
\put(0,1){\line(1,1){4}}
\qbezier(4,5)(7,8)(8,1)
\qbezier(4,5)(6.5,5)(8,3)
\put(0,1){\line(2,1){8}}
\put(0,5){\line(2,-1){8}}
\put(4,3){\line(1,0){4}}
\end{picture}$$
The three vertical edges correspond from the left respectively to $\{y\eq 0\}$, $\{x\eq 0\}$, $\{z\eq 0\}$. Indeed, we have $d\eq12$, and $|L^{[4]}|=9$; more precisely
$$L^{[4]}=L^{[(4)]}=\{x^3\pl y^3\pl z^3\eq xyz\eq 0\}\subset\PP^2,$$
and ${\rm Sing}\,L=L^{[2]}\sqcup L^{[4]}$, see a calculation in Example~(3.5) below. These imply that
$$E_{L,4}=\tfrac{d}{4}=3.$$
Here $H^1(\Ff,\C)_{\la}$ with $\ord\,\la\eq4$ is ESV-calculable, and its dimension is 2 (since it is a 4-net). The subset $J$ can be given by any factor in the last factorization of $f$.
It is known that $H^1(\Ff,\C)_{\la}\eq 0\,$ ($\la\,{\notin}\,\mu_4$), and $\dim H^1(\Ff,\C)_{\la}\eq 2\,$ ($\la\,{\in}\,\mu_4\1{\setminus}\1\{1\}$), see \cite{BDS}. (This can be confirmed by using a small computer program calculating the pole order spectral sequence.)
\msn
{\bf Example~3.5.} For $m\ges3$, let
$$f=\mprod_{i,j=0}^{m-1}\,(\zeta^ix\pl\zeta^jy\pl z),$$
with $\zeta:=e^{2\pi\sqrt{-1}/m}$. (This is a subarrangement of the Hessian arrangement when $m\eq 3$.)
We have $d\eq m^2$, $|L^{[m]}|=|L^{[(m)]}|=3\1m$ with $L^{[m]}=L\cap\{xyz\eq 0\}$, and ${\rm Sing}\,L=L^{[2]}\cup L^{[m]}$. Indeed, the last assertion can be reduced to the injectivity of the map
$$\Theta_i{\times}\Theta_j\ni(\la,\la')\mapsto\la/\la'\in\C,$$
with $\Theta_i:=\{\la\in\C\mid |\la\pl\zeta^i|\eq1,\,\la\ne 0\}$, calculating the intersection points of lines in $L$. (Here we get that $x=-(\zeta^{j'}\mi\zeta^j)/(\zeta^{i'}\mi\zeta^i)$ with $i,j$ fixed after setting $y\eq 1$.)
W may assume $i\eq j\eq0$ using the action of $\mu_m$ on $\C$. We have
$$\Theta_0=\bl\{re^{i\theta}\,\big|\,\,r=-2\cos\theta,\,\,\theta\in\bl(\tfrac{1}{2}\pi,\tfrac{3}{2}\pi\br)\br\}.$$
If $\la_1/\la'_1\eq\la_2/\la'_2$, we have $\arg\la_1\mi\arg\la'_1\eq\arg\la_2\mi\arg\la'_2$ (which is denoted by $\al\in(0,\pi)$).
Setting $\,\theta\eq\arg\la'$ (so that $\,\arg\la\eq\theta\pl\al$), the assertion is then reduced to the injectivity of the map
$$\bl(\tfrac{1}{2}\pi,\tfrac{3}{2}\pi\mi\al\br)\ni\theta\mapsto\cos(\theta\pl\al)/\cos\theta\in\R,$$
where the image is equal to $\,\cos\al\mi\sin\al\tan\theta$. So the assertion follows. We thus get that
$$\tfrac{d}{m}\eq m,\q E_{L,m}\eq 3.$$
\sk
In the $m$ even case, we can verify that $H^1(\Ff,\C)_{\la}$ for $\ord\,\la\eq m$ is {\it not\1} ESV-calculable as follows: Assume $J$ is given by
$$J=\bl\{([a_k],[b_k])\in(\Z/m\Z)^2\mid k\in\Z/m\1\Z\br\},$$
with $a_k,b_k\in[0,m{-}1]$. Here $(i,j)\in(\Z/m\Z)^2$ is identified with a line defined by the equation
$$\zeta^ix\pl\zeta^jy\pl z\eq0.$$ 
In order to satisfy condition~(a) or (b), we must have
$$\{a_k\}=\{b_k\}=\{a_k\mi b_k\pl m\delta_k\}=\{0,\dots,m{-}1\},$$
taking the intersection with $\{y\eq0\}$, $\{x\eq0\}$, $\{z\eq 0\}$ respectively.
Here $\delta_k=1$ if $a_k<b_k$, and 0 otherwise. We then get that
$$\msum_{k=0}^{m-1}\,a_k=\msum_{k=0}^{m-1}\,b_k=\msum_{k=0}^{m-1}\,(a_k\mi b_k\pl m\delta_k)=m(m{-}1)/2.$$
However, these imply that $\msum_k\,\delta_k=(m{-}1)/2$, which is a contradiction, since $m$ is even and $\delta_k\in\Z$. So ESV-non-calculability follows.
\sk
In the $m$ odd case, it is easy to see that $H^1(\Ff,\C)_{\la}$ for $\ord\,\la\eq m$ is ESV-calculable, since $J$ can be given by
$$J=\bl\{(i,2i)\in(\Z/m\Z)^2\mid i\in\Z/m\1\Z\br\}.$$
Note that the multiplication by 2 is an automorphism of $\Z/m\1\Z$ (with $m$ odd).
\msn
{\bf Example~3.6.} Consider a general hyperplane section of a reflection arrangement of type $G_{31}$. Here $d\eq 60$ and each line has 16 triple points and 3 points of multiplicity 6, see \cite{BDY}, \cite{OT}, \cite{ac}. Hence
$$E_{L,3}=19,\q E_{L,6}=3.$$
In this case $H^1(\Ff,\C)_{\la}$ with $\ord\,\la\eq3$ or 6 is never ESV-calculable, see {\it loc.\,cit.}

\end{document}